\documentclass[12pt]{article}
\usepackage{graphicx} 
\usepackage{lakshay}
\usepackage[backend=biber,style=alphabetic]{biblatex}
\usepackage[margin=1.25in]{geometry}
\usepackage{tikz}
\usepackage{quiver}
\usepackage{wrapfig}
\usepackage{matlab-prettifier}
\usepackage{indentfirst}
\usetikzlibrary{arrows}

\renewenvironment{abstract}
  {\small\quotation\noindent\textnormal{\MakeUppercase{\abstractname}. }\ignorespaces}
  {\endquotation}

\title{\large\bfseries\MakeUppercase{Evolution of Internal Waves on Surfaces of Revolution with Flat Ends}}
\author{\normalsize LAKSHAY PATEL}
\date{}

\begin{document}
\maketitle

\begin{abstract}
    We generalize the 2D linear internal wave equation to non-compact surfaces of revolution with flat ends. We show that the leading profile of this equation with oscillatory forcing is the unique incoming solution to the corresponding stationary equation. This is done by proving a limiting absorption principle via propagation of semiclassical defect measure.
\end{abstract}
\section{Introduction}

Let $\Omega$ be the surface of revolution \begin{equation}\label{eq:defofomeg}
    \Omega = \{(x_1, x_2, x_3)\in \R^3: x_1\in \R, x_2^2+x_3^2 = r_0^2(x_1)\} 
\end{equation}
where $r_0(x_1)\in C^\infty(\R; (0, \infty))$ and $r_0(x_1)$ is identically $1$ outside a compact set of $x_1\in \R$. Let $(x, \theta)\in \R\times \mathbb{S}^1$ be coordinates for the arc-length parametrization of $\Omega$ and $r(x) = r_0(x_1)$ the radius in terms of $x$. In these coordinates, the induced metric on $\Omega$ from $\R^3$ is 
\begin{equation}
    g = dx^2 + r(x)^2d\theta^2
\end{equation}
and the Laplace-Beltrami operator on $\Omega$ is
\begin{equation}
    \Delta_{\Omega} = \frac{1}{r(x)}\partial_{x}(r(x)\partial_{x}) + \frac{1}{r(x)^2}\partial_{\theta}^2
\end{equation}
We will henceforth use these coordinates to identify $\Omega$ with $\R\times \mathbb{S}^1$.

We consider the following Poincar\'e problem on $\Omega$:
\begin{equation}\label{eq:poincareomega}
    \begin{cases}
        \left(\partial_t^2\Delta_{\Omega} +
        r(x)^{-2}\partial_{\theta}^2\right)u = f(x, \theta)\cos(\lambda_0 t) \\
        u|_{t=0} = \partial_t u|_{t=0} = 0
    \end{cases}
\end{equation}
where $f\in C_c^\infty(\Omega)$ and $\lambda_0\in (0,1)$. Equation \eqref{eq:poincareomega} is a generalization of the 2D linear internal wave equation with oscillatory forcing, where $u$ is the stream function for a density-stratified Boussinesq fluid. For further discussion, see equation \eqref{eq:poincareog}. 

Our main result describes the long-time asymptotics of solutions to \eqref{eq:poincareomega}:

\begin{thm}\label{thm:longtimebehavior}
    Suppose $f\in C_c^\infty(\Omega)$ has no zero modes in $\theta$. Then the solution to \eqref{eq:poincareomega} is of the form
    \begin{equation}\label{eq:timeasymp}
        u(t) = Re(e^{i\lambda_0 t}u^+) + R(t) + \mathcal{E}(t)
    \end{equation}
    where $u^+\in H^{1, -\alpha}(\Omega)$ for all $\alpha>1/2$, $R(t)$ is uniformly bounded in $H^1(\Omega)$ with respect to $t$, and $\mathcal{E}(t)$ tends to zero in $H^{0,-\beta}(\Omega)$ as $t\to \infty$ for all $\beta >3/2$. Furthermore, the leading profile $u^+$ is the unique incoming solution to the stationary equation
    \begin{equation}
        (\lambda_0^2\Delta_\Omega + r(x)^{-2}\partial_{\theta}^2)u^+ = f
    \end{equation}
\end{thm}
\noindent Section \ref{sec:preliminaries} defines the Sobolev spaces used in theorem \ref{thm:longtimebehavior} and the notion of an incoming stationary solution.

We approach equation \ref{eq:poincareomega} as follows: 
the operator $\Delta_\Omega$ is invertible
when restricted to the Sobolev space $H^{-1}_0(\Omega)$, where
\begin{equation}
    H^s_0(\Omega) = \{u\in H^s(\Omega,d\text{Vol}_g): \hat{u}(x, 0)\equiv 0 \text{ for all } x\in \R\}
\end{equation}
and $\hat{u}$ is the partial Fourier series of $u$ in the $\theta$ variable. For $w(t) =\Delta_\Omega u(t)$, equation \eqref{eq:poincareomega} becomes
\begin{equation}\label{eq:evprob}
    \begin{cases}
        (\partial_t^2 + P)w = f(x, \theta)\cos(\lambda_0 t) \\
        w|_{t=0} = \partial_t w|_{t=0} = 0
    \end{cases}   
\end{equation}
where 
\begin{equation}
    P = r(x)^{-2}\partial_{\theta}^2\Delta_{\Omega}^{-1} : H^{-1}_0(\Omega) \to H^{-1}_0(\Omega)
\end{equation}
To solve equation \eqref{eq:evprob}, we make use of the functional calculus for $P$ by equipping $H^{-1}_0(\Omega)$ with the inner product
\begin{equation}\label{eq:h-10inprod}
    \langle u, v \rangle_{H^{-1}_0(\Omega)} = \langle \nabla_g \Delta_{\Omega}^{-1} u, \nabla_g\Delta_{\Omega}^{-1} v \rangle_{L^2(\Omega)}
\end{equation}
Then $P$ is a bounded self-adjoint operator on $H_0^{-1}(\Omega)$ with spectrum $[0,1]$. Using the functional calculus for $P$, $w(t)$ can be written in terms of the solution to the corresponding scalar equation and the spectral measure $dE(z)$ for $P$. To prove theorem \ref{thm:longtimebehavior}, we establish mapping properties and regularity of the spectral measure $dE(z)$. Stone's theorem expresses $dE(z)$ as the difference of limiting resolvents, which are
\begin{equation}\label{eq:reslimits}
    \lim_{\epsilon\downarrow 0}(P-\lambda^2 \pm i\epsilon)^{-1}
\end{equation}
Our core argument develops a \textit{limiting absorption principle}, meaning we show the limits \eqref{eq:reslimits} exist between suitable Sobolev spaces (theorem \ref{thm:lap}) and are differentiable in the spectral parameter $\lambda$ (theorem \ref{thm:diffinomg}).

\subsection{Relation to Internal Waves Literature.} Equation \eqref{eq:poincareomega} is a generalization of the 2D linear internal wave equation on a domain $\Omega \subset \R^2_{x_1, x_2}$, which is given by 
\begin{equation} \label{eq:poincareog}
\begin{cases}
    (\partial_t^2\Delta +\partial_{x_2}^2)u = f\cos(\lambda_0 t) \\
    u|_{t=0} = u_t|_{t=0} = 0 \\
    u|_{\partial \Omega} = 0
\end{cases}
\end{equation}
with $f\in C_c^\infty(\Omega)$ and $\lambda_0\in (0,1)$. Here $u(t)$ represents the stream function of linear inertial waves forced by an oscillatory profile (see for instance \cite{ralston-1973}). The formation of geometric patterns in internal waves was observed experimentally by Maas and Lam in \cite{maas-lam-1995}. Motivated by these observations, Colin de Verdi\`{e}re and Saint-Raymond used spectral properties of an associated zeroth order operator to explain the formation of such patterns \cite{colin-de-verdiere-saint-raymond-2019}. This was later improved by Dyatlov-Zworski in the context of radial point estimates \cite{Dyatlov_2019}. Dyatlov-Wang-Zworski then proved asymptotics for equation \eqref{eq:poincareog} in a 2D aquarium with smooth boundary \cite{dyatlov-wang-zworski-2025}.

In our case, the surface of revolution $\Omega$ generalizes an unbounded channel with flat ends, which is modeled by
\begin{equation}\label{eq:channelomega}
    \Omega_{\R^2} = \{(x_1, x_2) : x_1\in \R, G(x_1) - \pi \leq x_2 \leq 0\}
\end{equation}
where $G\in C_c^\infty(\R; (-\infty, \pi))$. Li-Wang-Wunsch investigate the scattering of internal waves on $\Omega_{\R^2}$ in \cite{li2025internalwaves2dsubcritical}, \cite{li2025evolutioninternalwaves2d}. In \cite{li2025evolutioninternalwaves2d}, they use reduction to the boundary to prove a limiting absorption principle. The surface of revolution $\Omega$ defined in \eqref{eq:defofomeg} is a boundaryless generalization of the unbounded channel with flat ends, where the angular variable $\theta$ replaces the vertical direction $x_2$ in \eqref{eq:channelomega}. The Sobolev space $H^{-1}_0(\Omega)$ and the assumption that $f$ has no zero modes in $\theta$ replace the Dirichlet boundary condition in \eqref{eq:poincareog}. The term $r(x)^{-2}\partial_{\theta}^2$ appearing in \eqref{eq:poincareomega} replaces the $\partial_{x_2}^2$ term in \eqref{eq:poincareog} as the angular part of the Laplace-Beltrami operator $\Delta_\Omega$. Our analysis of the boundaryless problem exhibits many common features of the methods in previous work, as well as a different method for obtaining a limiting absorption principle.

\subsection{Structure of the Paper.}
To show that $\lim_{\epsilon \downarrow 0} P(\lambda\pm i\epsilon)$ exists, we prove uniform estimates for the norm of $P(\lambda\pm i \epsilon)$ as an operator on suitable Sobolev spaces. First, these estimates are derived for the unperturbed/free operator ($g\equiv 1$) in section \ref{sec:flatcylinder}. The free estimates are then used in section \ref{sec:hifreq} and \ref{sec:lowfreq} to derive estimates for the perturbed operator in different $\theta$-frequency regimes. Our argument for high frequencies adapts Burq's argument in \cite{burq-2002} by using propagation of semiclassical defect measure to prove a resolvent estimate. In section \ref{sec:lap} we use these estimates and a compactness argument to obtain the desired limiting absorption principle and regularity of the limiting profile $u^+$. Then we show in section \ref{sec:diff} that the limiting resolvents are differentiable in $\lambda$ and use this in section \ref{sec:asymp} to prove theorem \ref{thm:longtimebehavior}.

\subsection*{Use of AI.} The author used a fusion of Copilot models to assist with spell check and LaTeX formatting. All mathematical ideas are either the author's own or from cited works. The author takes full responsibility for the mathematical correctness, proofs, and validity of this paper. 

\subsection*{Acknowledgments.} The author thanks Andy Badea, Hans Christianson, Perry Kleinz, Zhenhao Li, and especially Jared Wunsch for helpful discussions. The author was partially supported by NSF grant DMS-2136217.

\section{Preliminaries and Notation}\label{sec:preliminaries}

\subsection{Reduction of the Stationary Equation.}
Let $\omega\in \mathbb{C}$. For convenience we make a few modifications to the equation
\begin{equation}
    (P-\omega^2)w = f
\end{equation} 
We add $\Delta_\Omega$ back to the equation and set $u = \Delta_{\Omega}w$ to make this a second order equation:
\begin{equation}
    (r(x)^{-2}\partial_{\theta}^2 - \omega^2\Delta_{\Omega})u = f
\end{equation}
We exchange the first order term $\partial_{x}$ in $\Delta_\Omega$ for a potential by conjugating by $r(x)^{1/2}$:
\begin{equation}
    r^{1/2}(x)\Delta_\Omega r^{-1/2}(x) = \partial_{x}^2 + r^{-2}(x)\partial_{\theta}^2 + V(x)
\end{equation}
where 
\begin{equation}
    V(x) = \frac{1}{4}\left(\frac{r'(x)}{r(x)}\right)^2 - \frac{1}{2}\frac{r''(x)}{r(x)} \in C_c^\infty(\R_{x}, \R)
\end{equation}
is a smooth, compactly supported potential. By our assumptions on $r$, conjugating by $r^{1/2}$ is an isomorphism on any desired Sobolev space and fixes the $\partial_{\theta}^2$ term. So henceforth we work with the following equivalent equation:
\begin{equation}\label{eq:modifstationeqn}
    (\omega^2D_x^2 + g(x)(\omega^2-1)D_\theta^2 - \omega^2 V)u = f
\end{equation}
where we set
\begin{equation}
    g(x) = \frac{1}{r(x)^2}
\end{equation}
for convenience. From our assumptions on $r$, $g\in C^\infty(\R, (0, \infty))$ and $g(x)\equiv 1$ outside a compact set in $x$. For $\omega\in \mathbb{C}$, define 
\begin{equation}
    P(\omega) = \omega^2D_x^2 + g(x)(\omega^2-1)D_\theta^2 - \omega^2 V
\end{equation}
which is a family of second order operators and is hyperbolic when $\omega=\lambda\in (0,1)$. Then it suffices to develop a limiting absorption principle for $P(\omega)$ as $\omega\to \lambda$ from above or below the real axis.

\subsection{Incoming and Outgoing Solutions.} First, note that if $\Im(\omega) \neq 0$, then $P(\omega)$ is an elliptic operator. By a standard Lax-Milgram argument, $P(\omega): H_0^2(\Omega) \to L^2_0(\Omega)$ is invertible. Invertibility between these spaces fails as $\Im(\omega)\downarrow  0$, so we introduce the weighted Sobolev spaces
\begin{flalign}
    H^{s, \alpha}(\Omega) = \{u\in \mathcal{D}'(\Omega): \langle x\rangle^{\alpha}u \in H^s(\Omega)\}
\end{flalign}
Outside of the supports of $V, g-1,$ and $f$, solutions to $P(\lambda)u=f$ must have the form
\begin{equation}
    \hat{u}(x, k) = \begin{cases}
            a_-(\lambda,k)e^{i|k|c_\lambda x}  + b_-(\lambda,k)e^{-i|k|c_\lambda x} & x \ll 0 \\
            a_+(\lambda,k)e^{-i|k|c_\lambda x} + b_+(\lambda,k)e^{i|k|c_\lambda x}  & x \gg 0
        \end{cases}
\end{equation}
where 
\begin{equation}
         c_\lambda = \sqrt{\lambda^{-2}-1} 
\end{equation}
and $a_\pm(\lambda, k), b_\pm(\lambda, k)$ are constants. To distinguish between such solutions, we make the following definition:

\begin{dfn}
    Let $\alpha > 1/2$. A solution $u\in H^{0, -\alpha}(\Omega)$ of $P(\lambda)u = f\in L^2_{comp}(\Omega)$ is called \textit{incoming} if each $\theta$-Fourier mode $\hat{u}(x, k)$ has the form
    \begin{equation}
        \hat{u}(x, k) = \begin{cases}
            a_-(\lambda,k)e^{i|k|c_\lambda x} & x \ll 0 \\
            a_+(\lambda,k)e^{-i|k|c_\lambda x}& x \gg 0
        \end{cases}
    \end{equation}
    for constants $a_-(\lambda, k), a_+(\lambda, k)$. A solution $u\in H^{0, -\alpha}(\Omega)$ of $P(\lambda)u = f$ is called \textit{outgoing} if each $\theta$-Fourier mode $\hat{u}(x, k)$ has the form
    \begin{equation}
        \hat{u}(x, k) = \begin{cases}
            b_-(\lambda,k)e^{-i|k|c_\lambda x} & x \ll 0 \\
            b_+(\lambda,k)e^{i|k|c_\lambda x} & x \gg 0
        \end{cases}
    \end{equation}
    for constants $b_-(\lambda, k), b_+(\lambda, k)$.
\end{dfn}

Note that the partial Fourier coefficients $\hat{u}(x, k)$ solve
\begin{equation}\label{eq:1dhelm}
    (D_x^2 - c_\lambda^2g(x)k^2 - V(x))\hat{u}(x, k) = \lambda^{-2}\hat{f}(x, k)
\end{equation}
Outside a compact set, \eqref{eq:1dhelm} is a 1D Helmholtz equation. The definitions above agree with the standard incoming/outgoing conditions on $\hat{u}(x, k)$ for the 1D Helmholtz equation.

\begin{lem} \label{lem:uniqincsol}
    Let $f\in L^2_{comp}(\Omega)$ have no zero modes in $\theta$. Then there is at most one incoming (or outgoing) solution $u\in \mathcal{D}'(\Omega)$ to $P(\lambda)u = f$ such that $u$ has no zero modes in $\theta$.
\end{lem}

\begin{pf}
    Let $u$ be the difference between any two incoming solutions for the same source term $f$, which is incoming and satisfies the homogeneous equation $P(\lambda)u = 0$. For each $k$ the $k$-th Fourier mode $\hat{u}(x, k)$ solves the homogeneous 1D equation
    \begin{equation}
        (D_x^2 - c_\lambda^2g(x)k^2 - V(x))\hat{u}(x, k) = 0
    \end{equation}
    Let $v_k(x) = \hat{u}(x, k)$. Let $J(x) =  \Im(\bar{v}_k\partial_{x}v_k)$ be the Wronskian. Then 
    \begin{align}
        2iJ'(x) &= \bar{v}_k\partial_{x}^2v_k - v_k\partial_{x}^2\bar{v}_k \\
        &= \bar{v}_k(c_\lambda^2g(x)k^2 + V(x))v_k - v_k(c_\lambda^2g(x)k^2 + V(x))\bar{v}_k \\
        &= 0
    \end{align}
    so $J(x)$ is constant. Since $v_k$ is an incoming solution to the 1D equation, there are constants $a_-, a_+$ such that 
    \begin{equation}
        v_k(x) = \begin{cases}
            a_-e^{i|k|c_{\lambda}x} & x \ll 0 \\
            a_+e^{-i|k|c_\lambda x} & x \gg 0
        \end{cases}
    \end{equation}
    Then $J(x)=|k|c_{\lambda}|a_-|^2$ for $x\ll 0$ and $J(x)=-|k|c_{\lambda}|a_+|^2$ for $x\gg 0$. Since $J(x)$ is constant, $a_-=a_+=0$. Since $v_k$ satisfies the 1D equation and is compactly supported, $v_k\equiv 0$ and therefore $u\equiv 0$ if $\hat{u}(x, 0) \equiv 0$. The same argument applies if $u$ is the difference of two outgoing solutions.
    \qed
\end{pf}

\begin{rem}
    In sections \ref{sec:flatcylinder} - \ref{sec:diff} we will only state and prove results for $\lambda + i\epsilon$ with $\epsilon \downarrow 0$. This corresponds to incoming solutions to the stationary equation, but results in these sections also hold for the case $\epsilon \uparrow 0$ and correspond to outgoing solutions.
\end{rem}
    
\section{Flat Cylinder} \label{sec:flatcylinder}

We begin by studying the flat cylinder ($r\equiv 1$). In this case $g\equiv 1$ and $V\equiv 0$, so the operator in equation \ref{eq:modifstationeqn} becomes the constant coefficient operator
\begin{equation}
    P_0(\omega) = \omega^2D_x^2 + (\omega^2-1)D_\theta^2
\end{equation}
In this section, we prove the following theorem:
\begin{thm}\label{flatres}
    Suppose $f\in L^2_{comp}(\Omega)$ and $\hat{f}(x, 0)\equiv 0$. Let $u(\omega)$ solve 
    \begin{equation}
        P_0(\omega) u(\omega) = f
    \end{equation}
    Then for any $\alpha > 1/2$, $u(\omega)$ converges in $H^{0, -\alpha}(\Omega)$ to the unique incoming solution $u(\lambda)$ of $P_0(\lambda)u(\lambda)  =f$ as $\omega \downarrow \lambda$. Further, for any $K\subset (0,1)$ compact there is $C>0$ depending only on $K$ such that
    \begin{equation}\label{eq:freeP0est}
        \|u(\omega)\|_{H^{0, -\alpha}} \leq C \|f\|_{L^2}
    \end{equation}
    for $\Re(\omega)\in K, Im(\omega)\geq 0$. 
\end{thm}

As in the proof of lemma \ref{lem:uniqincsol}, we note that the partial Fourier series coefficients of $u(\omega)$ in $\theta$ satisfy the equation
\begin{equation}
    (\omega^2D_x^2 + (\omega^2-1)k^2)\hat{u}(\omega,x, k) = \hat{f}(x, k)
\end{equation}
To prove both this case and the perturbed case, we create a semiclassical 1-D problem. Set $h=|k|^{-1}$ and divide by $k^2\omega^{-2}$, which yields
\begin{equation}
    (h^2D_x^2 -(\omega^{-2}-1))\hat{u}(\omega,h) = h^2 \omega^{-2}\hat{f}(h)
\end{equation}

\begin{rem}
    We write $\hat{u}(\omega, h)$ and $\hat{f}(h)$ even though there are two frequencies $\pm k$ corresponding to $h$. In what follows we take $k>0$ and the case $k<0$ is identical.
\end{rem}

We will first study the limit of the operator $(h^2D_x^2 - (\omega^{-2}-1))^{-1}$ as $\epsilon \downarrow 0$ before summing Fourier coefficients to prove theorem \ref{flatres}. Let
\begin{equation}
    Q_0(\omega, h) = h^2D_x^2 - (\omega^{-2}-1)
\end{equation}
and
\begin{equation}
    R_0(\omega,h) = Q_0(\omega, h)^{-1} : L^2_0(\R) \to H^2(\R)
\end{equation}
Note that $(\omega^{-2}-1)$ converges to $\lambda^{-2}-1$ as $\epsilon\downarrow 0$ and $\Im(\omega^{-2}-1) <0$. We set 
\begin{equation}
    c_{\omega} = \sqrt{\omega^{-2}-1}
\end{equation}
to be the branch with $\Im(c_\omega)>0$, so that
\begin{equation}
    \lim_{\omega \downarrow \lambda}c_{\omega} = -\sqrt{\lambda^{-2}-1} = -c_\lambda
\end{equation}
We recall standard results for the 1D semiclassical Helmholtz equation and include proofs for completeness. 

\begin{lem}\label{lem:1dsemlap}
    Let $w\in L^2_{comp}(\R)$ and set $v(\omega, h) := R_0(\omega, h)w$. Then $v(\omega, h)$ converges in $H^{0, -\alpha}(\R)$ to the unique incoming solution of
    \begin{equation}
        (h^2D_x^2 - c_\lambda^2)v(\lambda, h) = w
    \end{equation}
\end{lem}

\begin{rem}
    We will generally use $v(\omega, h)$ to denote solutions to the semiclassical 1D equation in sections \ref{sec:flatcylinder} - \ref{sec:diff}, with the application $v(\omega, h) = \hat{u}(\omega,k)$ in mind. 
\end{rem}

\begin{pf}
    The $L^2$ integral kernel of $R_0(\omega, h)$ is 
    \begin{equation}\label{eq:freereskernelnonweighted}
        R_0(\omega, h; x, y) = \frac{i}{2hc_{\omega}}e^{\frac{i}{h}c_{\omega}|x-y|}
    \end{equation}
    Let $\chi\in C_c^{\infty}(\R)$ be a smooth cutoff function on the support of $w$. Then set 
    \[ \tilde{R_0}(\omega, h; x, y) = \langle x\rangle^{-\alpha}R_0(\omega, h; x, y)\chi(y)\]
    For fixed $h$, $R_0(\omega, h; x, y)$ converges in $L^2(\R_x\times \R_y)$ to
    \begin{equation}\label{eqn:freeweightkern}
        \tilde{R_0}(\lambda, h; x, y) = -\langle x\rangle^{-\alpha}\frac{i}{2hc_\lambda}e^{-\frac{i}{h}c_\lambda|x-y|}\chi(y)
    \end{equation}
    as $\omega \downarrow \lambda$. The corresponding operators $\tilde{R_0}(\omega, h)$ converge to $\tilde{R_0}(\lambda, h)$ in Hilbert-Schmidt norm and therefore in operator norm $L^2\to L^2$. This means that $v(\omega, h)$ converges to
    \begin{equation}
        v(\lambda, h) := \langle x\rangle^{\alpha}\tilde{R_0}(\lambda, h)w
    \end{equation}
    in $H^{0, -\alpha}(\R)$. Take $R>0$ large enough so that $w$ is supported in $(-R, R)$. Then $v(\lambda, h)$ is of the form 
    \begin{equation}
        v(\lambda, x, h) = 
        \begin{cases}
            a_-e^{\frac{i}{h}c_\lambda x} & x\ll 0 \\
            a_+ e^{-\frac{i}{h}c_\lambda x}& x\gg 0
        \end{cases}
    \end{equation}
    where 
    \begin{equation}
        a_- = \int_{-\infty}^{R} e^{-\frac{i}{h}c_\lambda y_1} w(y_1) \, dy_1, \quad a_+ = \int_{-R}^{\infty} e^{\frac{i}{h}c_\lambda y_1} w(y_1) \, dy_1
    \end{equation}
    So $v(\lambda, h)$ is the unique incoming solution.
    \qed
\end{pf}

To prove the estimate in theorem \ref{flatres} and establish $H^1$ regularity of $u$, we improve the convergence in lemma \ref{lem:1dsemlap} to $H^{1, -\alpha}$ and keep track of the $h$ dependence:

\begin{lem}\label{lem:freeresest}
    Let $K\subset (0,1)$ be compact. Then for $\lambda\in K$, $R_0(\omega, h)w$ converges in $H^{1, -\alpha}(\R)$ to $R_0(\lambda, h)w$ as $\omega\downarrow \lambda$ and is bounded by
    \begin{equation}\label{eq:freeresh1alph}
        \|R_0(\lambda, h)w\|_{H^{1, -\alpha}_h} \leq \frac{C}{h}\|w\|_{L^2}
    \end{equation}
    where $C$ is independent of $h,\lambda$ and the norm on the left-hand side is the semiclassical Sobolev norm
    \begin{equation}
        \|v\|_{H^{1, -\alpha}_h} = \|\langle x\rangle^{-\alpha}v\|_{L^2} + \|hD(\langle x\rangle^{-\alpha}v)\|_{L^2}
    \end{equation}
\end{lem}

\begin{pf}
    We bound the kernel \eqref{eqn:freeweightkern} by 
    \begin{equation}
        \|\tilde{R_0}(\lambda, h; x, y)\|_{L^2(\R_x\times \R_y)}^2 \leq \frac{C}{hc_\lambda}
    \end{equation}
    For $K\subset (0,1)$ compact, $c_\lambda$ is bounded away from $0$, so $\|\tilde{R_0}(\lambda, h)\|_{L^2\to L^2} \leq Ch^{-1}$ for some $C$ independent of $\lambda$. This gives the $L^2$ part of the estimate. For the $hD$ part, we differentiate the kernel which yields
    \begin{equation}
    \begin{aligned}
        hD_x\tilde{R_0}(\omega, h; x, y) &= -\frac{i}{2c_\omega}\langle x\rangle^{-\alpha}\text{sgn}(x-y)e^{-\frac{i}{h}c_\omega|x-y|}\chi(y) \\ 
        &- \frac{i\alpha}{2hc_\omega}\langle x\rangle^{-\alpha-1}e^{-\frac{i}{h}c_\omega|x-y|}\chi(y)
    \end{aligned}
    \end{equation}
    The right-hand side converges in $L^2(\R_x\times \R_y)$ to the corresponding kernel for $hD\tilde{R_0}(\lambda, h; x, y)$, and
    \begin{equation}
        \|hD_x\tilde{R_0}(\omega, h; x, y)\|_{L^2(\R_x\times \R_y)}^2 \leq \frac{C}{hc_\omega}
    \end{equation}
    So convergence holds in $H^{1, -\alpha}(\R)$ with estimate \eqref{eq:freeresh1alph}. \qed
\end{pf}

We now prove theorem \ref{flatres}:

\begin{pf}
    Write $\hat{u}(\omega, x, k)$ in terms of the free resolvent as:
    \begin{equation}
        \hat{u}(\omega, k) = R_0(\omega, h)\left(h^2\omega^{-2}\hat{f}(-, k)\right)
    \end{equation}
    By Lemma \ref{lem:freeresest} each Fourier mode $\hat{u}(\omega)(x,k)$ converges to some $v(\lambda, x, k)$ in $H^1(\R)$ as $\omega\downarrow \lambda$. $v$ satisfies the estimates
   \begin{equation}\label{freel2est}
       \|v(\lambda, -,k)\|_{H^{0,-\alpha}(\R)} \leq \frac{C}{|k|} \| \hat{f}(-,k) \|_{L^2(\R)}
   \end{equation}
   \begin{equation}\label{freeh1est}
       \|D_{x}v(\lambda, -,k)\|_{H^{0,-\alpha}(\R)} \leq C \| \hat{f}(-,k) \|_{H^{0, \alpha}(\R)}
   \end{equation}
   Since the norms $\|\hat{f}(k)\|_{H^{0, \alpha}(\R)}$ are square-summable and we have estimate $\ref{freel2est}$, it follows that
   \begin{equation}
    \hat{u}(\omega,x,k) \to v(\lambda, x, k) \quad \text{in } \langle x\rangle ^{\alpha}L^2(\R_{x}\times \Z_k) \quad \text{as } \omega \downarrow \lambda.
   \end{equation}
   Set $u(\lambda)$ to be the Fourier sum of the coefficients $v(\lambda, k)$. Then estimate \ref{freel2est} implies that 
   \begin{equation}
    \|\langle x\rangle^{-\alpha}D_\theta u(\lambda) \|_{L^2} \leq  C \| f\|_{L^2}
   \end{equation}
   Estimate \ref{freeh1est} implies that 
   \begin{equation}
    \|\langle x\rangle^{-\alpha}D_x u(\lambda) \|_{L^2} \leq  C \| f\|_{L^2}
   \end{equation}
   So the desired $H^{1, -\alpha}$ estimate follows. For each Fourier mode $k$, $\hat{u}(\lambda, k) = v(\lambda,k)$ is the incoming solution to
   \begin{equation}
    (\lambda^2 D_x^2 +(\lambda^2-1)k^2)\hat{u}(\lambda, x, k) = \hat{f}(x, k)
   \end{equation}
   So $u(\lambda)$ is the unique incoming solution to
   \begin{equation}
    P_0(\lambda)u(\lambda) = f
   \end{equation}
\qed
\end{pf}

\section{High Frequency Estimate} \label{sec:hifreq}

We now take $r(x)$ to be a general function satisfying $r>0$ and $r-1\in C_c^{\infty}(\R_{x})$. As in section \ref{sec:flatcylinder} we set $h=|k|^{-1}$. Then the partial Fourier coefficients $\hat{u}(\omega, x, k)$ satisfy the 1D semiclassical equation
\begin{equation}\label{eq:1dperturbedeqn}
    (h^2D_x^2 - g(x)c_\omega^2 - h^2V(x))\hat{u}(\omega,h) = h^2\omega^{-2}\hat{f}(h)
\end{equation}
Let
\begin{equation}
    Q(\omega, h) := h^2D_x^2 - g(x)c_\omega^2 - h^2V(x)
\end{equation}
and
\begin{equation}
    R(\omega, h) := Q(\omega, h)^{-1}
\end{equation}
Equation \ref{eq:1dperturbedeqn} is similar to the 1D semiclassical Helmholtz equation, except that the perturbation $c_\omega^2(1-g)(x)$ is weighted by the spectral parameter $c_\omega$. We still refer to the operator $R(\omega, h)$ as a \textit{resolvent} to emphasize connections to standard resolvent theory.

In this case we do not have an explicit kernel for $R(\omega, h)$. So we obtain convergence of $\hat{u}(\omega,h)$ by proving an $L^2_{loc}$ analogue of the uniform estimate in lemma \ref{lem:freeresest}. In this section we prove this estimate for high frequencies (small $h$):

\begin{thm}\label{thm:hifreqest}
    Let $\chi \equiv 1$ on the support of $g-1, V$ and $K\subset (0,1)$ compact. Then for all $N\in \mathbb{N}$ there is $h_N>0$ and $C>0$ such that for $0<h<h_N$, $0 < \Im(\omega) \leq Nh$ and $\Re(\omega)\in K$ we have the cutoff resolvent estimate
    \begin{equation}\label{eq:hifreqest}
        \|\chi R(\omega, h)\chi\|_{L^2 \to L^2} \leq \frac{C}{h}
    \end{equation}
    where $C$ only depends on $K, N$.
\end{thm}
\begin{rem}
    Unlike the case of the standard semiclassical Helmholtz equation, theorem \ref{thm:hifreqest} is true regardless of the critical points of $g$. This is due to the fact that the perturbation is weighted by the spectral parameter $c_\omega$, see lemma \ref{lem:suppmu} and the discussion following lemma \ref{lem:invar}.
\end{rem}

We prove theorem \ref{thm:hifreqest} by propagation of semiclassical defect measure, which is an adaptation of Burq's argument in theorem 1.3 of \cite{burq-2002}. The proof is by contradiction: assuming estimate \ref{eq:hifreqest} is false for some N, there are sequences $h_n\to 0$, $\omega_n$ with $0\leq \Im(\omega_n) \leq Nh_n$, and $f_n\in L^2$ such that 
\[ \|\chi R(\omega_n, h_n)\chi f_n\|_{L^2} \geq \frac{n}{h_n}\|f_n\|_{L^2}\]
Set $u_n = R(\omega_n, h_n)(\chi f_n)$
and normalize $\|\chi u_n\|_{L^2}$ so that
\begin{equation}\label{eq:defmeassetup}
    \begin{cases}
        \|\chi u_n\|_{L^2} = 1, \\
        \|\chi f_n\|_{L^2} \leq \frac{h_n}{n} = o(h_n)
    \end{cases}
\end{equation}
Since $R(\omega_n, h_n)$ is invertible, we know a priori that $u_n\in L^2$ although no uniform bound on $\|u_n\|_{L^2}$. Before constructing a semiclassical defect measure for $\{u_n\}$ we first show that $\{u_n\}$ is bounded in $L^2_{loc}$:
\begin{lem}\label{lem:solnl2loc}
    For any $\tilde{\chi} \in C_c^{\infty}(\R)$, $\|\tilde{\chi}u_n\|_{L^2}$ is uniformly bounded. 
\end{lem}

\begin{pf}
    Let $\psi\in C^\infty(\R)$ be such that $\psi\equiv 0$ on the supports of $g-1$ and $V$, $\chi\equiv 1$ on the support of $D\psi$, and $\psi\equiv 1$ for large $|x|$. Since $\|\chi u_n\|=1$ for all $n$, it suffices to show that $\psi u_n$ is uniformly bounded in $L^2_{loc}$. For this we will use the free operator $Q_0$ and resolvent $R_0(\omega_n, h_n)$. Since $Q\equiv Q_0$ on the support of $\psi$, 
    \begin{flalign}
        Q_0(\omega_n, h_n)\psi u_n &= Q(\omega_n, h_n)\psi u_n \\
        & = \psi \chi f_n + [Q(h_n, \omega_n), \psi]u_n \\
        & = \psi \chi f_n + h_n^2D\psi Du_n - h_n^2\psi'' u_n \label{eq:rhsofpsiun}
    \end{flalign}
    Since lemma \ref{lem:freeresest} bounds the norm of $R_0(\omega, h): L^2_{comp} \to L^2_{loc}$ by $Ch^{-1}$, it suffices to show that the right-hand side \ref{eq:rhsofpsiun} is bounded by $O(h_n)$ in $L^2_{loc}$. We already have this for $\psi \chi f_n$ and $h_n^2 \psi '' u_n$ from \ref{eq:defmeassetup}. For the $h_nDu_n$ term, first note that $D\psi Du_n = \chi D\psi Du_n$. We integrate by parts to get
    \begin{flalign}
        \|\chi D\psi h_n Du_n\|^2 &= \langle \chi^2(D\psi)^2h_n^2Du_n, Du_n\rangle \\
        &= \langle 2\chi D\chi(D\psi)^2 h_n^2 Du_n, u_n\rangle \label{eq:ibp11} \\
        &+ \langle 2\chi^2 (D^2\psi)D\psi h_n^2 Du_n, u_n\rangle \label{eq:ibp12} \\
        &+ \langle \chi^2(D\psi)^2h_n^2 D^2u_n, u_n \rangle \label{eq:ibp13}
    \end{flalign}
    The first term \ref{eq:ibp11} is zero because $D\psi$, $D\chi$ have disjoint supports. For the third term \ref{eq:ibp13} we use the fact that $h_n^2D^2u_n = c_{\omega_n}^2 u_n + \chi f_n$:
    \begin{flalign}
        |\langle \chi^2(D\psi)^2h_n^2 D^2u_n, u_n \rangle| &\leq |\langle \chi^2(D\psi)^2 (c_{\omega_n}^2u_n+h_n^2V), u_n\rangle| + |\langle \chi^2(D\psi)^2 \chi f_n, u_n\rangle| \\
        &\leq |c_{\omega_n}|^2\|\chi u_n\|^2 + \|\chi f_n\|\|\chi u_n\| \\
        &= O(1) + o(h_n) \\
        & = O(1)
    \end{flalign}
    We estimate the second term \ref{eq:ibp12} by the Peter-Paul inequality:
    \begin{flalign}
        |\langle \chi^2 (D^2\psi)(D\psi)h_n^2D u_n, u_n\rangle| &= h_n|\langle \chi (D\psi)h_n Du_n,(D^2\psi)\chi u_n\rangle| \\
        &\leq h_n \left(\frac{\delta}{2}\|\chi(D\chi)h_n Du_n\|^2 + \frac{1}{2\delta}\|(D^2\psi)\chi u_n\|^2\right) \\
        &\lesssim h_n\delta\|\chi(D\chi)h_n Du_n\|^2 + O(h_n)
    \end{flalign}
    Taking $\delta>0$ small enough to absorb $\delta\|\chi h_n Du_n\|$ into the left-hand side and putting these estimates together shows that $\|\chi D\psi h_n Du_n\|$ is $O(1)$. It follows that
    \begin{equation}
        \|\psi \chi f_n + h_n^2D\psi Du_n - h_n^2\psi'' u_n\|_{L^2} = O(h_n)
    \end{equation}
    Now take any $\tilde{\chi}\in C_c^\infty$ supported in the support of $\psi$. The weighted free resolvent bound implies that $\|\tilde{\chi}\psi u_n\|$ is $O(1)$, which finishes. \qed
\end{pf}

We now construct a semiclassical defect measure associated to $u_n$.

\begin{lem}\label{lem:defmeas}
    There is a subsequence $\{u_{n_k}\}$ and a positive Radon measure $\mu$ on $T^*\R$ such that for any $a\in C_c^\infty(T^*\R)$, 
    \begin{equation}
        \lim_{k\to\infty} \langle a(x, h_{n_k}D)u_{n_k}, u_{n_k}\rangle = \int ad\mu
    \end{equation}
    where $a(x, hD)$ is the semiclassical left-quantization of $a$:
    \begin{equation}
        a(x, hD)u(x) = \frac{1}{2\pi h}\int\int e^{\frac{i}{h}(x-y)\xi}a(x, \xi)u(y)dyd\xi
    \end{equation}
\end{lem} 

See \S~5.1 of \cite{zworski-2012} for the construction of $\mu$. We extract further subsequences so that
\begin{equation}
    h_{n_k}^{-1}\Im(\omega_{n_k}) \to \alpha\in [0, N]
\end{equation}
\begin{equation}
    \omega_{n_k} \to \lambda \in (0,1)
\end{equation}
 For convenience, index these sequences by $n$ instead of $n_k$. We now show that the measure $\mu$ satisfies certain properties from which we will deduce a contradiction, beginning with the support of $\mu$. 

\begin{lem}\label{lem:suppmu}
    The support of $\mu$ is contained in the set 
    \begin{equation}
        \Sigma = \{(x, \xi) : \xi^2 = c_{\lambda}^2g(x)\}
    \end{equation}
    which is the characteristic set of the semiclassical principal symbol of $Q(\lambda, h)$:
    \begin{equation}
        q(x, \xi) = \xi^2 - c_\lambda^2g(x)
    \end{equation}
\end{lem}

\begin{pf}
    Let $a\in C_c^\infty(T^*\R)$. Then 
    \begin{equation}
        \langle a(x, h_nD)Q(\lambda, h_n)u_n, u_n\rangle = \langle a(x, h_nD)\chi f_n, u_n\rangle + (c_\lambda^2 - c_{\omega_n}^2)\langle a(x, h_nD)u_n, gu_n\rangle 
    \end{equation}
    Since $u_n$ is bounded in $L^2_{loc}$ and $\|\chi f_n\| = o(h_n)$, both terms on the right-hand side tend to zero. The operator $a(x, h_nD)Q(\lambda, h_n)$ has semiclassical principal symbol $a(x, \xi)q(x, \xi)$. By the semiclassical symbol calculus
    \[ a(x, hD)Q(\lambda, h) = (aq)(x, hD) + O_{L^2_{comp}\to L^2_{loc}}(h) \]
    So $\langle \mu, aq\rangle = 0$ for all $a$, which implies that $\mu$ is supported in the set $\{q =0\} = \Sigma$. \qed
\end{pf}

The conditions on $u_n, f_n$ in \ref{eq:defmeassetup} imply that $\mu$ is nontrivial:
\begin{lem} \label{hnosc}
    Let $\phi(\xi)\in C_c^\infty(\R_\xi)$ be a smooth cutoff function such that $\phi(\xi)\equiv 1$ on the fibre of $\Sigma$ over the support of $\chi$. Then
    \begin{equation}
        \mu(\chi(x)^2\phi(\xi)^2) = 1
    \end{equation}
\end{lem}

\begin{pf} 
    Let $a(x,\xi) = \chi(x)^2\phi(\xi)^2$ and $b(x, \xi) = \chi(x)\phi(\xi)$. The semiclassical principal symbol of $a(x, hD)$ is the same as that of $b(x, hD)^*b(x, hD)$, so 
    \begin{align}
        \mu(\chi^2(x)^2\phi(\xi)^2) &= \lim_{n\to \infty} \langle a(x, h_nD) u_n, u_n\rangle \\
        &= \lim_{n\to \infty}\|b(x, h_nD)u_n\|^2
    \end{align}
    For the same reason the right quantization $b(h_nD, y)$ gives the same limit as $h_n\to 0$, so 
    \begin{align}
        \mu(\chi^2(x)^2\phi(\xi)^2) &= \lim_{n\to \infty}\|b(h_nD, y)u_n\|^2 \\
        &= \lim_{n\to \infty}\|\phi(\xi)\widehat{\chi u}(\xi)\|_{L^2}^2
    \end{align}
    By Lemma \ref{lem:suppmu} this is independent of the choice of $\phi$ as long as $\phi\equiv 1$ on the fibre over the support of $\chi$. Using the normalization of $\chi u_n$ in \eqref{eq:defmeassetup}, we let the support of $\phi$ exhaust $\R_{\xi}$ and find that 
    \begin{align}
        \mu(\chi^2(x)^2\phi(\xi)^2) &= \lim_{n\to \infty}\|\phi(\xi)\widehat{\chi u}(\xi)\|_{L^2}^2 \\
        &= \lim_{n\to \infty}\|\chi u_n\|^2 \\
        & = 1
    \end{align} \qed
\end{pf}

We eventually show that $\mu$ must be zero by using the free estimates from section \ref{sec:flatcylinder}. For this we first show that $\mu$ is invariant under the Hamiltonian flow of the symbol $q(x, \xi)$:
\begin{lem}\label{lem:invar}
    For $a\in C_c^\infty(T^*\R)$, we have 
    \[\langle \mu, H_qa\rangle = 4\alpha c_\lambda\langle\mu,g\cdot a\rangle\]
    where 
    \[ H_q(a) = \partial_\xi q(x, \xi)\partial_x a(x, \xi) - \partial_x q(x, \xi)\partial_\xi a(x, \xi) \]
    In particular, the support of $\mu$ is invariant under the flow of $H_q$
\end{lem}

\begin{pf}
    Let $a\in C_c^\infty(T^*\R)$. By the semiclassical symbolic calculus,
    \begin{equation}
        ih_n^{-1}[a(x, h_nD), Q(\lambda, h_n)] = H_qa(x, h_nD) + O_{L^2_{loc}\to L^2_{comp}}(h_n)
    \end{equation}
    To make use of the equation $Q(h_n, \omega_n)u_n = \chi f_n$, we rewrite the commutator as
    \begin{equation}
        [a(x, h_nD), Q(\lambda, h_n)] = [a(x, h_nD), Q(\omega_n, h_n)]  + [a(x, h_nD), (c_{\omega_n}^2 - c_{\lambda}^2)g(x)]
    \end{equation}
    Since $a,g$ are zeroth order semiclassical symbols, 
    \begin{equation}
        [a(x, h_nD), g(x)] = O_{L^2_{loc}\to L^2_{comp}}(h_n)
    \end{equation}
    Since $c_{\omega_n}^2 - c_\lambda^2 = o(1)$, the second term 
    \begin{equation}
        ih_n^{-1}[a(x, h_nD), (c_{\omega_n}^2 - c_{\lambda}^2)g(x)]
    \end{equation}
    vanishes as $n\to\infty$. So it suffices to calculate the limit of the first term:
    \begin{flalign}
        ih_n^{-1}\langle [a(x, h_nD), Q(\omega_n, h_n)]u_n, u_n\rangle &=ih_n^{-1}\langle a(x, h_nD)\chi f_n, u_n\rangle \\ &- ih_n^{-1}\langle a(x, h_nD)u_n, Q^*(\omega_n, h_n)u_n\rangle \\
        &= ih_n^{-1}\langle a(x, h_nD)\chi f_n, u_n\rangle \\ &- ih_n^{-1}\langle a(x, h_nD)u_n, \chi f_n\rangle \\ &- ih_n^{-1}(c_{\omega_n}^2 - \overline{c_{\omega_n}}^2)\langle a(x, h_nD)u_n, gu_n\rangle 
    \end{flalign}
    Since $\|\chi f_n\| = o(h_n)$, the first two terms converge to $0$. We rewrite the last term as
    \begin{flalign}
        ih_n^{-1}(c_{\omega_n}^2 - \overline{c_{\omega_n}}^2)\langle a(x, h_nD)u_n, gu_n\rangle &= -2h_n^{-1}\Im(c_{\omega_n}^2)\langle g(x)a(x, h_nD)u_n, u_n\rangle \\
    \end{flalign}
    Recall that $\inv{h_n}\Im(c_{\omega_n}) \to \alpha$, so 
    \begin{flalign}
        \lim_{n\to\infty} ih_n^{-1}(c_{\omega_n}^2 - \overline{c_{\omega_n}}^2)\langle a(x, h_nD)u_n, gu_n\rangle &= 4\alpha c_\lambda\langle \mu, g\cdot a\rangle
    \end{flalign} \qed
\end{pf}

We now investigate the dynamics of the Hamiltonian system given by $q$. The Hamilton vector field of $q$ is
\begin{equation}
    H_q = 2\xi\partial_x + 2c_\lambda^2g'(x)\partial_\xi
\end{equation}
So the characteristic set $\Sigma$ has two connected components,
\begin{equation}
    \Sigma_{\pm} = \left\{\left(x, \pm c_\lambda\sqrt{g(x)}\right)\right\}
\end{equation}
On each component of $\Sigma$ the sign of $\xi$ remains fixed, so the sign of $H_q(x)$ is fixed along the $H_q$ flow. Therefore the flow of every point in $\Sigma$ escapes any compact subset of $x$ in both the forward and backward time direction. In particular, every point in the fibre of $\Sigma$ over the support of $\chi$ will eventually reach the free resolvent region under the $H_q$ flow. 

To make this more precise, let $R>0$ be large enough so that $\text{supp}(\chi) \subset \{|x| < R\}$. Then the forward time flow of each point in the fibre over the support of $\chi$ eventually reaches the outgoing set, defined by
\begin{equation}
    Out = \{(x, \xi)\in \Sigma : |x|>R, x\cdot \xi > 0\}
\end{equation}
We use the properties of the free resolvent $R_0(\omega_n, h_n)$ to show the following:
\begin{lem}\label{lem:outnull}
    $\mu(Out) = 0$
\end{lem}
\begin{pf}
    Suppose $a\in C_c^\infty(T^*\R)$ is supported in $Out$. Take $\psi$ to be the same function as in Lemma \ref{lem:solnl2loc}. We showed $\psi u_n = R_0(\omega_n, h_n)\tilde{f_n}$ with $f_n$ bounded in $L^2$. So
    \begin{equation}
        \psi u_n = \begin{cases}
            a_-(\omega_n, h_n)e^{-ic_{\omega_n}x/h_n} & x < -R \\
            a_+(\omega_n, h_n)e^{ic_{\omega_n}x/h_n} & x > R
        \end{cases}
    \end{equation}
    where $R>0$ is sufficiently large. Note that although this looks like an outgoing solution, $c_{\omega_n}\to -c_\lambda$ by our convention. Since $\psi\equiv 1$ on the base projection of the support of $a$, 
    \begin{equation}
        \lim_{n\to \infty}\langle a(x, h_nD)u_n, u_n\rangle = \lim_{n\to\infty}\langle a(x, h_nD)\psi u_n, \psi u_n\rangle
    \end{equation}
    One can then check using stationary phase that the right-hand side converges to $0$, which implies that $\mu(a) = 0$. \qed
\end{pf}
\newline

\noindent The proof of theorem \ref{thm:hifreqest} now follows from putting these lemmas together. By lemma \ref{lem:invar}, the support of $\mu$ is invariant under the $H_q$ flow. By lemma \ref{lem:outnull}, $\mu$ has no mass on the outgoing set. Since every point in the fibre over the support of $\chi$ flows to the outgoing set, this implies that $\mu(\chi(x)^2\phi(\xi)^2) = 0$, which contradicts lemma \ref{hnosc}.

\section{Low Frequency Estimate} \label{sec:lowfreq}

The estimate in theorem \ref{thm:hifreqest} is only valid for sufficiently small $h$. For low frequencies $h>h_N$ we obtain the same estimate via a perturbative argument. Since there are finitely many $h > h_N$, it suffices to prove the estimate for fixed $h$. If $\eta = c_\omega/h$ then
\begin{align}
    (h^2D^2 - c_\omega^2g(x) -h^2V(x))^{-1} &= h^{-2}(D^2 - (c_\omega/h)^2g(x) - V(x))^{-1} \\
    &= h^{-2}(D^2 - \eta^2g(x) - V(x)) 
\end{align}
so we consider the following operators on $\R$:
\begin{align}
    Q(\eta) &:= D^2 - \eta^2g(x) - V(x)\\
    R(\eta) &:= Q(\eta)^{-1}
\end{align}
Using analytic Fredholm theory, we will show that $R(\eta)$ is a meromorphic family of operators $L^2_{comp}(\R)\to L^2_{loc}(\R)$ for $\Im(\eta)\leq 0$ with no poles on the real axis. Regularity of $R(\eta)$ will be used in section \ref{sec:diff} as well. Let
\begin{equation}
    R_0(\eta) = (D^2 - \eta^2)^{-1}
\end{equation}
$R_0(\eta)$ is well-defined for $\Im(\eta)<0$. By inspection of the Schwartz kernel for $R_0(\eta)$ (which can be read off from equation \ref{eq:freereskernelnonweighted}), $R_0(\eta) : L^2_{comp}\to L^2_{loc}$ extends meromorphically to $\mathbb{C}$ with a single pole at $\eta = 0$. We will treat $Q(\eta)$ as a compact perturbation of $Q_0(\eta)$ to prove the following:

\begin{lem}\label{lem:resismero}
    For any $\chi\in C_c^\infty(\R)$, $\chi R(\eta)\chi: L^2(\R) \to L^2(\R)$ is well-defined for $\Im(\eta) < 0$ and extends to be meromorphic for $\eta\in \mathbb{C}$. 
\end{lem}

\begin{pf}
    Set $\tilde{g}(x) = g(x)-1$. Then 
    \begin{equation}\label{eq:compactperturb}
        Q(\eta) = (D^2 - \eta^2) + \eta^2\tilde{g}(x) + V(x)
    \end{equation}
    Since $\tilde{g}$ is compactly supported, $\eta^2\tilde{g}(x) + V(x) : H^2(\R)\to L^2(\R)$ is a compact operator. Therefore $Q(\eta)$ is a Fredholm operator of index $0$ for $\Im(\eta)<0$. Since $Q(\eta)$ has no kernel for $\Im(\eta)<0$, it is invertible. Applying the operator $R_0(\eta)$ to equation \ref{eq:compactperturb} yields 
    \begin{equation}
        Q(\eta)R_0(\eta) = I + (\eta^2\tilde{g}(x)+V(x))R_0(\eta)
    \end{equation}
    Let 
    \begin{equation}
        K(\eta) = (\eta^2\tilde{g}(x)+V(x))R_0(\eta)
    \end{equation}
    Then $\chi K(\eta)\chi : L^2 \to L^2$ is a meromorphic family of operators for $\eta\in \mathbb{C}\setminus\{0\}$. Without loss of generality $\chi\equiv 1$ on the support of $\tilde{g}, V$, so
    \begin{equation}
        (I + K(\eta)(1-\chi))^{-1} = I - K(\eta)(1-\chi)
    \end{equation}
    and
    \begin{flalign}
        (I + K(\eta))^{-1} &= (I + K(\eta)(1-\chi) + K(\eta)\chi)^{-1} \\
        &= \left[(I + K(\eta)(1-\chi))(I+K(\eta)\chi)\right]^{-1} \\
        &= (I+K(\eta)\chi)^{-1}(I-K(\eta)(1-\chi)) \label{eq:fred1}
    \end{flalign}
    Since $K(\eta)\chi$ is compact, $(I+K(\eta)\chi)^{-1}$ is meromorphic in $\eta$ by the analytic Fredholm theorem. So $(I+K(\eta))^{-1}$ is meromorphic in $\eta$, and therefore $\chi R(\eta)\chi = \chi R_0(\eta)\chi(I+K(\eta))^{-1}$ is meromorphic in $\eta$. \qed
\end{pf}

Next we show that $\chi R(\eta)\chi$ has no poles on the real axis.
\begin{lem}\label{lem:norealpoles}
    Each $\eta\in \R_{>0}$ is a regular point of $\chi R(\eta)\chi$ and gives the incoming solution to $Q(\eta)u = f$. 
\end{lem}
\begin{rem}
    $\eta\in\R_{>0}$ corresponds to $\lambda\in (0,1)$ in the original spectral parameter. 
\end{rem}

\begin{pf}
Fix $f\in L^2_{comp}(\R)$ and set $u(\eta) = R(\eta)f$. For $\Im(\eta)<0$, 
\begin{align}
    u(\eta) &= R(\eta)f \\
    &= R_0(\eta)(I+K(\eta))^{-1}f \\
    &= R_0(\eta)\tilde{f}
\end{align}
where $\tilde{f}\in L^2_{comp}(\R)$. Since $Q_0(\eta)$ gives an incoming solution to the free equation, $u(\eta)$ is an incoming solution to $Q(\eta)u = f$. Let $\psi_{\pm}\in C^\infty(\R)$ such that $\psi\equiv 0$ on the support of $\tilde{g}, f$, $\psi_+\equiv 1$ for $x\gg 0$, and $\psi_-\equiv 1$ for $x\ll 0$. Then the incoming condition is equivalent to
\begin{equation}
    \psi_{\pm}(D \pm \eta)u(\eta) =0
\end{equation}
This condition extends to all regular $\eta$ by meromorphic continuation. From the analytic Fredholm theorem, $(I + K(\eta))^{-1}$ has a pole at $\eta$ if and only if $I + K(\eta)$ has a nontrivial kernel. If $u\in \ker(I + K(\eta))$, then $u$ will be a nontrivial incoming solution to $Q(\eta)u = 0$. But no such solution exists, so $\eta$ is a regular point of $(I+K(\eta))^{-1}$ and therefore of $\chi R(\eta)\chi$. \qed
\end{pf}
 
We can now improve theorem \ref{thm:hifreqest} to all frequencies:

\begin{thm} \label{thm:unifresest}
    Let $K\subset (0,1)$ be compact and $\chi \equiv 1$ on the support of $g-1$. Then for all $N\in \mathbb{N}$ there is $C>0$ depending only on $K,N$ such that for $\Im(\omega) \leq Nh$ and $\Re(\omega)\in K$ we have the cutoff resolvent estimate
    \begin{equation}\label{eq:unifresest}
        \|\chi(h^2D^2 - c_{\omega}^2g(x) - h^2V)^{-1}\chi\|_{L^2 \to L^2} \leq \frac{C}{h}
    \end{equation}
\end{thm}

\begin{pf}
    From theorem \ref{thm:hifreqest} there is $h_N>0$ such that the estimate holds for $0<h<h_N$. For $h>h_N$, we have 
    \begin{equation}
        \|\chi(h^2D^2 - c_{\omega}^2g(x) - h^2V(x))^{-1}\chi\| = h^{-2}\|\chi R(c_\omega/h)\chi\|
    \end{equation}
    For this range of $\omega$ and $h$, the spectral parameter $\eta = c_\omega/h$ lies in a compact neighborhood where $\chi R(\eta)\chi$ is holomorphic. So $\|\chi R(c_\omega/h)\chi\|$ is bounded by some constant $C'>0$, which implies
    \begin{equation}
        h^{-2}\|\chi R(c_\omega /h)\chi\| \leq \frac{C'}{h^2} \leq \frac{C'}{h_N}\cdot \frac{1}{h}
    \end{equation}
    for $h_N < h<1$. So for the constant $\max{C, C'/h_N}$ the estimate \ref{eq:unifresest} holds for all $h>0$. \qed
\end{pf}

\section{Limiting Absorption Principle} \label{sec:lap}

Before extracting a limit of $u(\omega)$, we make some improvements to the estimate in theorem \ref{thm:unifresest}. For this section we use $v(h, x)$ to denote solutions to the 1D semiclassical equation, and eventually set $v(h, x) = \hat{u}(x, k)$ for the limiting absorption principle of equation \ref{eq:modifstationeqn}. We begin by improving the regularity of estimate \ref{eq:unifresest} to $H^1$:
\begin{cor}\label{cor:h1hboost}
    Let $H^1_h(\R)$ be the semiclassical Sobolev space with norm
    \begin{equation}
        \|v(h)\|_{H^1_h} = \|v(h)\|_{L^2} + \|hD(v(h))\|_{L^2}    
    \end{equation}
    Let $K\subset (0,1)$ be compact. Suppose $v(\omega, h)$ is a family of solutions to 
    \begin{equation}
        (h^2D^2 - c_\omega^2 g - h^2V)v(\omega, h) = f
    \end{equation} 
    with $f\in L^2_{comp}$. Let $\chi\in C_c^\infty(\R)$ with $\chi\equiv 1$ on the support of $g-1, f, V$. Then for all $N>0$ there is $C>0$ such that for $\Im(\omega) \leq Nh$ and $\Re(\omega)\in K$ the solution $v(\omega, h)$ satisfies the estimate
    \begin{equation}\label{eq:h1hboost}
        \|\chi v(\omega)\|_{H^1_h} \leq \frac{C}{h}\|f\|_{L^2}
    \end{equation} 
\end{cor}

\begin{pf}
    For convenience we drop the argument $\omega$ from $v$. It suffices to estimate $\|hD(\chi v)\|_{L^2}$. By the product rule,
    \begin{equation}
        \|h D(\chi v)\|^2 \leq \|h(D\chi)v\|^2 + \|h\chi (Dv)\|^2
    \end{equation}
    Let $\tilde{\chi}\in C_c^\infty(\R)$ be such that $\tilde{\chi}\equiv 1$ on the support of $\chi$. The first term is controlled by $h\|\tilde{\chi}v\|$, which is controlled in $L^2_{loc}$ by estimate \ref{eq:unifresest}. For the $Dv$ term we integrate by parts and use the Peter-Paul inequality to find
    \begin{flalign}
        \|h \chi (Dv)\|^2 &= \int h^2\chi^2Dv\cdot Dv \\
        &= 2\int(\chi hDv)(h(D\chi) v) + \int (\chi h^2D^2 v)(\chi v) \\
        &\leq 4\delta \|h\chi Dv\|^2 + \frac{h}{\delta}\|\tilde{\chi}v\| + \left(\|\chi f\| + c_\omega^2\|g\chi v\| + h^2\|\chi V v\|\right)\|\chi v\| \\
        &\lesssim \delta \|h\chi Dv\|^2 + \frac{1}{h^2}\|f\|^2 \label{eq:rhs00}
    \end{flalign}
    By taking $\delta>0$ small enough we absorb the $\|h\chi Dv\|$ term of \ref{eq:rhs00} into the left-hand side and find that
    \begin{equation}
        \|h\chi (Dv)\|^2 \leq\frac{C}{h}\|f\|^2
    \end{equation}
    \qed
\end{pf}

Next we improve this to a global weighted $H^1$ estimate:
\begin{cor}\label{cor:loctoweigh}
    Let $\alpha>1/2$ and $K\subset (0,1)$ compact. Suppose $v(\omega, h)$ solves 
    \begin{equation}
        (h^2D^2 - c_\omega^2 g - h^2V)v(\omega, h) = f
    \end{equation}
    with $f\in L^2_{comp}$. Then for all $N>0$ there is $C>0$ such that for $\Im(\omega) \leq Nh$ and $\Re(\omega)\in K$ the solution $v(\omega, h)$ satisfies the estimate
    \begin{equation}\label{eq:wtest}
       \|v(\omega)\|_{H^{1,-\alpha}_h} \leq \frac{C}{h}\|f\|_{L^2} 
    \end{equation}
\end{cor}

\begin{pf}
    Let $\psi = 1-\chi$. Since we have the local estimate \ref{eq:h1hboost} it suffices to prove \ref{eq:wtest} for $\psi v$. On the support of $\psi$ the operator $h^2D^2-c_\omega^2g - h^2V$ equals the free operator $h^2D^2 - c_\omega^2$. So 
    \begin{flalign}
        (h^2D^2 - c_\omega^2g - h^2V)(\psi v) &= (h^2D^2 - c_\omega^2g)(\psi v) \\
        &= f + (h^2D^2\psi)\cdot v + h^2(D\psi)(Dv)
    \end{flalign}
    From the estimates in the proof of corollary \ref{cor:h1hboost} and lemma \ref{lem:solnl2loc}, this right-hand side is controlled by $\|f\|_{L^2}$ uniformly in $h, \omega$. Then
    \begin{equation}
        \psi v = R_0(\omega, h)\left(f + (h^2D^2\psi)\cdot v + h^2(D\psi)(Dv)\right)
    \end{equation}
    From section \ref{sec:flatcylinder} the norm of  $R_0(\omega, h): L^2_{comp} \to H^{1, -\alpha}_h$ is controlled by $Ch^{-1}$, so estimate \ref{eq:wtest} follows. \qed
\end{pf}

We now have all the pieces to prove the following theorem:
\begin{thm}\label{thm:lap}
    Let $f\in L^2_{comp}(\Omega)$ and $\lambda\in (0,1)$. For $\epsilon>0$ let $u(\lambda+i\epsilon)\in L^2(\Omega)$ be the unique solution to 
    \begin{equation}
        P(\lambda+i\epsilon)u({\lambda+i\epsilon}) = f
    \end{equation}
    Then for any $\alpha>1/2$, $u(\lambda+i\epsilon) \to u(\lambda)$ in $H^{1, -\alpha}(\Omega)$ as $\epsilon\downarrow 0$, where $u(\lambda)$ is the unique incoming solution to $P(\lambda)u(\lambda) = f$. For any compact subset $K\subset (0,1)$ there is $C>0$ depending only on $K$ such that 
    \begin{equation}
    \|u(\lambda)\|_{H^{1,-\alpha}(\Omega)} \leq C \|f\|_{L^2}
    \end{equation}
\end{thm}

\begin{pf}
    Fix a nonzero Fourier mode $k$ in $\theta$ and set $h=|k|^{-1}$. Then $\hat{u}(\omega, k)$ solves the 1D semiclassical equation 
    \begin{equation}
        Q(\omega, h)\hat{u}(\omega, h) = h^2\omega^{-2}\hat{f}(h)
    \end{equation}
    Choose some $\beta$ such that $1/2 <\beta<\alpha$ and corollary \ref{cor:loctoweigh} implies
    \begin{equation}
        \|\hat{u}(\omega, h)\|_{H^{1, -\beta}_h} \leq Ch\|\hat{f}(h)\|_{L^2}
    \end{equation}
    So for fixed $h$ this estimate implies $\hat{u}(\lambda+i\epsilon, h)$ is uniformly bounded in $H^{1, -\beta}_h$ for $\epsilon \leq Nh$. Now fix any subsequence $\epsilon_n\to 0$. Since $H^{1, -\beta}(\R)$ is compactly embedded in $H^{0, -\alpha}(\R)$, there is a subsequence $\hat{u}(\lambda+\epsilon_{n_j}, h)$ that converges in $H^{0, -\alpha}_h$ to some $v(\lambda, h)$. In particular $\hat{u}(\lambda+\epsilon_{n_j}, h)\to v(\lambda, h)$ in distribution, so $v$ solves
    \begin{equation}
        Q(\lambda, h)v(\lambda, h) = h^2\lambda^{-2}\hat{f}
    \end{equation}
    Since $v$ is obtained via the limiting incoming resolvent, $v$ is an incoming solution. Since the incoming solution is unique, the limit $v$ is independent of the subsequence $\epsilon_{n_j}$. Therefore $\hat{u}(\lambda+i\epsilon, k)$ converges to a unique incoming solution $v(\lambda, x, k)$ in $H^{0, -\alpha}$ as $\epsilon\downarrow 0$. The limit will satisfy the estimate
    \begin{equation}\label{eq:x2decay}
        \|v(\lambda,-, k)\|_{H^{0, -\alpha}} \leq \frac{C}{|k|}\|\hat{f}(-, h)\|_{L^2}
    \end{equation}

    Next we show convergence in $H^{1, -\alpha}$. The arguments used in the proof of lemma \ref{lem:solnl2loc} and corollary \ref{cor:h1hboost} give the estimate 
    \begin{equation}\label{eq:showedh1}
        \|v(\lambda, -, k)\|_{H^{1, -\alpha}} \leq C\|\hat{f}(-, k)\|_{L^2}
    \end{equation}
    So we now have boundedness in $H^{1, -\alpha}$ and finish by checking convergence in $H^{1, -\alpha}$. Set 
    \begin{equation}
        w_\epsilon(x) = \hat{u}(\omega,x, k) - v(\lambda,x, k)
    \end{equation}
    Then 
    \begin{equation}
        w_\epsilon \to 0, D_x^2w_{\epsilon}\to 0 \quad \text{in } H^{0, -\alpha}_h
    \end{equation}
    From \ref{eq:showedh1} we integrate by parts to find
    \begin{equation}
        \|Dw_{\epsilon}\|^2_{H^{0, -\alpha}(\R)} \lesssim \|D^2 w_{\epsilon}\|^2_{H^{0, -\alpha}(\R)}\|w_{\epsilon}\|_{H^{0, -\alpha}(\R)} +  \|w_{\epsilon}\|_{H^{0, -\alpha}(\R)}
    \end{equation}
    So $\|D w_{\epsilon}\|_{H^{0, -\alpha}(\R)} \to 0$ as $\epsilon\to 0$, which gives convergence in $H^{1, -\alpha}$.

    Now we have that $\hat{u}(\lambda+i\epsilon,-, k)$ converges to $v(x, k)$ in $H^{1, -\alpha}$ pointwise in $k$. From corollary \ref{cor:loctoweigh}, 
    \begin{equation}
        \|\hat{u}(\lambda+i\epsilon) - v(\lambda)\|_{\ell^2(\Z_k, H^{1, -\alpha}(\R_{x}))} \to 0 \quad \text{as } \epsilon \to 0
    \end{equation}
    Set $u(\lambda)$ to be the Fourier sum of the coefficients $v(x, k)$. Then $u(\lambda)$ is an incoming solution to $P_\lambda u_\lambda = f$, and $u(\lambda+i\epsilon) \to u(\lambda)$ in $L^2(S_{\theta}, H^{1, -\alpha}(\R_{x}))$ as $\epsilon\to 0$. From estimate \ref{eq:x2decay} we also have $D_\theta\langle x\rangle^{-\alpha}u(\lambda+i\epsilon) \to D_\theta\langle x\rangle^{-\alpha}u(\lambda)$ in $L^2(\Omega)$. Combining these two gives the desired convergence in $H^{1, -\alpha}(\Omega)$ as well as the estimate.
    \qed
\end{pf}

\section{Differentiability in $\lambda$} \label{sec:diff}

In this section we show regularity of the incoming solution $u(\lambda)$ in $\lambda$. We use this in section \ref{sec:asymp} to show that the remainder term $\mathcal{E}(t)$ in theorem \ref{thm:longtimebehavior} tends to zero. 

\begin{thm}\label{thm:diffinomg}
    For all $\beta > 3/2$ the solution map $u(\lambda) : (0,1)\to H^{0, -\beta}(\Omega)$ is differentiable in $\lambda$. For $K\subset (0,1)$ compact, there is $C>0$ depending only on $K$ such that
    \begin{equation}\label{eq:diffestonomg}
        \|(\partial_\lambda u)(\lambda) \|_{H^{0, -\beta}(\Omega)} \leq C\|f\|_{L^2(\Omega)}
    \end{equation}
    for all $\lambda\in K$.
\end{thm}

As before, set $v(\lambda,h, x) = \hat{u}(\lambda, x, k)$. Since $v(\lambda, h, x)$ is incoming, it takes the form
\begin{equation} \label{eq:vlmbatends}
    v(\lambda,h, x) = \begin{cases}
    a_-(\lambda, h)e^{ic_\lambda x/h} & x \ll 0 \\
    a_+(\lambda, h)e^{-ic_\lambda x/h} & x \gg 0
    \end{cases}
\end{equation}
Differentiating these terms in $\lambda$ gives an extra factor of $\frac{x}{h}$, which explains the higher growth order $\beta > 3/2$ and the loss of differential order relative to theorem \ref{thm:lap}. To take advantage of the form \eqref{eq:vlmbatends} we first show that the free resolvent is differentiable:

\begin{lem}\label{lem:freeresdiff}
    For fixed $h$, the limiting incoming resolvent $R_0(\lambda, h) = (h^2D^2 - c_{\lambda}^2)^{-1}$ is differentiable as a map $(0,1)\to \mathcal{L}(L^2_{comp}(\R), H^{0, -\beta}(\R))$. 
\end{lem}

\begin{pf}
    Fix $\chi \in C_c^\infty(\R)$. Recall that the Schwartz kernel for $\langle x\rangle^{-\beta}R_0(\lambda, h)\chi(x)$ is 
    
    \begin{equation}
        R_0^{\beta, \chi}(\lambda; x, y) = -\frac{i}{2c_{\lambda}h}\langle x\rangle^{-\beta}e^{-ic_\lambda |x-y|/h}\chi(y)
    \end{equation}
    The kernels $R_0^{\beta, \chi}(\lambda)$ are differentiable in $\lambda$ from $(0,1)\to L^2(\R_x\times \R_y)$. So $\langle x\rangle^{-\beta}R_0(\lambda, h)\chi(x)$ is differentiable in Hilbert-Schmidt norm, hence operator norm $L^2\to L^2$. \qed
\end{pf}
\newline

For the perturbed case we first show that $v(\lambda,h)$ is differentiable for fixed $h$. Then we prove an estimate for $\|\partial_\lambda v(\lambda, h)\|_{H^{0, -\beta}}$ that is uniform in $h$, which we use to establish differentiability in the full space $H^{0, -\beta}(\Omega)$.

\begin{lem}\label{lem:difffixedh}
    For fixed $h$, the map $v(\lambda,h,x) : (0,1)\to H^{0, -\beta}(\R)$ is differentiable in $\lambda$. 
\end{lem}

\begin{pf}
    Apply the free operator $Q_0(\lambda, h)$ to $v(\lambda,h)$:
    \begin{equation}\label{eq:diffffixedh1}
        Q_0(\lambda, h)v(\lambda,h) = \lambda^{-2}h^2\hat{f}(h) + c_\lambda^2(g-1)v(\lambda,h) + h^2Vv(\lambda, h)
    \end{equation}
    Let $\chi\in C_c^\infty(\R)$ be a cutoff function that is identically 1 on the supports of $g-1, \hat{f}(h), V$. Then
    \begin{equation}
        c_\lambda^2(1-g)v(\lambda,h) = c_\lambda^2(1-g)\chi R(\lambda, h)\chi(\lambda^{-2}h^2\hat{f}(h))
    \end{equation}
    \begin{equation}
        h^2Vv(\lambda, h) = h^2V\chi R(\lambda, h)\chi(\lambda^{-2}h^2\hat{f}(h))
    \end{equation}
    By lemma \ref{lem:norealpoles}, the operator $\chi R(\lambda, h)\chi$ is holomorphic for $\lambda\in (0,1)$. Thus the right-hand side of equation \ref{eq:diffffixedh1} is differentiable $(0,1)\to L^2$ and is supported where $\chi\equiv 1$ for all $\lambda$. Applying $R_0(\lambda, h)$ to both sides of \ref{eq:diffffixedh1} and using lemma \ref{lem:freeresdiff} shows that $v(\lambda): (0,1)\to H^{0, -\beta}$ is differentiable. \qed
\end{pf}

Now we estimate the $h$-dependence of the norm $\|\partial_\lambda v(\lambda, h)\|_{H^{0, -\beta}}$:

\begin{lem}\label{lem:drvest}
    Let $K\subset (0,1)$ be compact and $\beta >3/2$. Then there is $C>0$ independent of $\lambda, h$ such that 
    \begin{equation}\label{eq:drvest}
        \|(\partial_\lambda v)(\lambda, h)\|_{H^{0, -\beta}} \leq C\|\hat{f}(h)\|_{L^2}
    \end{equation}
    for $\lambda \in K$. 
\end{lem}

\begin{pf}
    Our main strategy is to show that $\partial_\lambda v(\lambda, h)$ satisfies an equation involving $Q(\lambda,h)$ and use previous estimates to control its norm. We differentiate 
    \begin{equation}
        Q(\lambda, h)v(\lambda, h) = \lambda^{-2}h^2\hat{f}(h)
    \end{equation}
    with respect to $\lambda$ and find that $\partial_\lambda v$ solves 
    \begin{equation}\label{eq:diffedeqn}
        Q(\lambda, h)(\partial_\lambda v(\lambda,h)) = -2\lambda^{-3}g\cdot v(\lambda,h) - 2h^2 \lambda^{-3}\hat{f}(h)
    \end{equation}
    The right-hand side of this equation is not compactly supported owing to $g\cdot v(\lambda, h)$, so we will solve away the first term. We first show that the coefficients $a_\pm(\lambda, h)$ are differentiable and controlled in $h$, then differentiate the right-hand side of \ref{eq:vlmbatends}.
    
    Fix $R>0$ large enough so \ref{eq:vlmbatends} is valid for $\pm x\geq R$. Then
    \begin{equation}\label{eq:apmfromv}
        \int_R^\infty \langle x\rangle^{-2\beta}|v(\lambda, h, x)|^2dx = |a_+(\lambda, h)|^2\int_R^\infty\langle x\rangle^{-2\beta}dx
    \end{equation}
    Since $\langle x\rangle^{-\beta}v(\lambda, h, x)$ is differentiable in $L^2$, $a_+(\lambda, h)$ is differentiable as well. Equation \ref{eq:apmfromv} together with the estimate on $\|v(\lambda, h, \cdot)\|_{H^{0, -\beta}}$ from theorem \ref{thm:unifresest} also implies
    \begin{equation} \label{eq:coeffest}
        |a_+(\lambda, h)|\leq Ch\|\hat{f}(h)\|_{L^2}
    \end{equation}
    for $\lambda\in K$. The same argument holds for $a_-(\lambda, h)$, which justifies differentiating equation \ref{eq:vlmbatends}:
    \begin{equation}
        (\partial_\lambda v)(\lambda, h, x) = \begin{cases}
            a_-(\lambda, h)\frac{ix}{h}(\partial_\lambda c_\lambda) e^{ic_\lambda x/h} + \partial_\lambda a_-(\lambda, h)e^{ic_\lambda x/h} & x\ll 0 \\
            -a_+(\lambda, h)\frac{ix}{h}(\partial_\lambda c_\lambda) e^{-ic_\lambda x/h} + \partial_\lambda a_+(\lambda, h)e^{-ic_\lambda x/h} & x\gg 0 \\
        \end{cases}
    \end{equation}
    From estimate \ref{eq:coeffest} the leading order growth terms are controlled in $H^{0, -\beta}$ by $h\|\hat{f}(h)\|$. Fix $R>0$ large and choose $\psi_+\in C^\infty(\R)$ such that $\psi_+\equiv 1$ on $[R+1, \infty)$ and $\text{supp}(\psi_+)\subset (R, \infty)$. Set $\psi_-(x) = \psi_+(-x)$ and let
    \begin{equation}
        (\widetilde{\partial_\lambda v})(\lambda, h, x) = \psi_-(x)a_-(\lambda, h)\frac{ix}{h}(\partial_\lambda c_\lambda) e^{ic_\lambda x/h} - \psi_+(x)a_+(\lambda, h)\frac{ix}{h}(\partial_\lambda c_\lambda) e^{-ic_\lambda x/h}
    \end{equation} 
    for which we know 
    \begin{equation}\label{eq:tildederest}
        \|(\widetilde{\partial_\lambda v})(\lambda, h, x)\|_{H^{0, -\beta}} \leq C\|\hat{f}(h)\|_{L^2}.
    \end{equation}
    Let 
    \begin{equation}
        r(\lambda,h) = (\partial_\lambda v)(\lambda, h, \cdot) - (\widetilde{\partial_\lambda v})(\lambda, h, \cdot)
    \end{equation}
    Then $r(\lambda, h)$ solves the equation 
    \begin{equation}\label{eq:rlambdaeq}
        \begin{split}
            Q(\lambda, h)r(\lambda, h) &= (\psi_+ + \psi_- -1)2\lambda^{-3}gv(\lambda,h) - 2h^2\lambda^{-3}\hat{f}(h) \\
        &+[Q(\lambda, h), \psi_-]a_-(\lambda)(\partial_\lambda c_\lambda)\frac{ix}{h}e^{ic_\lambda x/h} \\
        &- [Q(\lambda, h), \psi_+]a_+(\lambda)(\partial_\lambda c_\lambda)\frac{ix}{h}e^{-ic_\lambda x/h} 
        \end{split}
    \end{equation}
    By construction,    
    \begin{equation}
        r(\lambda, h, x) = \begin{cases}
            \partial_\lambda a_-(\lambda, h)e^{ic_\lambda x/h} & x\ll 0 \\
            \partial_\lambda a_+(\lambda, h)e^{-ic_\lambda x/h} & x\gg 0
        \end{cases}
    \end{equation}
    so $r_\lambda(\lambda, h, x)$ is the unique incoming solution to equation \ref{eq:rlambdaeq}. The right-hand side of \ref{eq:rlambdaeq} is supported in $(-R-1, R+1)$ for all $\lambda, h$. By theorem \ref{thm:unifresest} the first term is bounded in $L^2$ by $h\|\hat{f}(h)\|$. The commutator appearing in the latter two terms is
    \begin{equation}
        [Q(\lambda, h), \psi_\pm] = h^2(D^2\psi_{\pm})+h(D\psi)\cdot (hD)
    \end{equation}
    so these terms are also bounded in $L^2$ by $h\|\hat{f}(h)\|$. Apply the resolvent estimate in theorem \ref{thm:unifresest} to the right-hand side of equation \ref{eq:rlambdaeq} and we obtain
    \begin{equation}
        \|r_\lambda(\lambda, h, \cdot)\|_{H^{0, -\beta}} \leq C\|f(h)\|_{L^2}
    \end{equation}
    for $\lambda\in K$. Combining this with estimate \ref{eq:tildederest} gives the desired result. \qed
\end{pf}
\newline

Returning to the full space $\Omega$, we can now prove theorem \ref{thm:diffinomg}. 
\newline

\begin{pf}
    Fix $\lambda\in (0,1)$ and consider the difference quotient 
    \begin{equation}
        \frac{\hat{u}(\lambda,x, k) - \hat{u}(\mu,x, k)}{\lambda - \mu}
    \end{equation}
    as $\mu\to \lambda$. Let $I$ be a compact interval in $(0,1)$ containing $\lambda$. Then, by the mean value inequality and lemma \ref{lem:drvest},
    \begin{equation}
        \left\|\frac{\hat{u}(\lambda,-,k) - \hat{u}(\mu,-,k)}{\lambda - \mu}\right\|_{H^{0, -\beta}(\R)} \leq \sup_{\mu \in I}\|(\partial_\lambda \hat{u})(\mu,-,k)\|_{H^{0, -\beta}(\R)} \leq C\|\hat{f}(-, k)\|_{L^2(\R)}
    \end{equation}
    where $C$ is independent of $k, \lambda\in K$. Since $\|\hat{f}(x,k)\|$ is square-summable in $k$ and these difference quotients converge pointwise for each $k$, the difference quotients converge in $\ell^2(\Z_k, H^{0, -\beta}(\R_{x}))$ to $(\partial_\lambda \hat{u})(\lambda,-,k)$. Summing the Fourier coefficients $(\partial_\lambda \hat{u})(\lambda,-,k)$ shows that $u_\lambda$ is differentiable in $H^{0, -\beta}(\Omega)$ and $\partial_\lambda u$ satisfies estimate \ref{eq:diffestonomg}. \qed
\end{pf}

\section{Long Time Behavior}\label{sec:asymp}
We now return to the time-dependent problem 
\begin{equation}
    \begin{cases}
        \left(\partial_t^2\Delta_{\Omega} +
        r(x)^{-2}\partial_{\theta}^2\right)u = f(x, \theta)\cos(\lambda_0 t) \\
        u|_{t=0} = \partial_t u|_{t=0} = 0
    \end{cases}
\end{equation}
and prove theorem \ref{thm:longtimebehavior}.

\begin{pf}
Recall the associated evolution equation
\begin{equation}
    \begin{cases}
        \left(\partial_t^2 +
        P\right)w = f(x, \theta)\cos(\lambda_0 t) \\
        w|_{t=0} = \partial_t w|_{t=0} = 0
    \end{cases}
\end{equation}
where
\begin{equation}
    P = r(x)^{-2}\partial_{\theta}^2\Delta_{\Omega}^{-1}
\end{equation}
One can show that for the inner product \eqref{eq:h-10inprod} on $H^{-1}_0$, $P$ is self-adjoint with spectrum $[0,1]$. We use the functional calculus to write $w(t)$ as
\begin{equation}
    w(t) = \Re(e^{i\lambda_0 t}W_{t,z}(P)f)
\end{equation}
where 
\begin{equation}\label{eq:wtdef}
    W_{t,\lambda_0}(z) := \int_0^t \frac{\sin(s\sqrt{z})}{\sqrt{z}}e^{-i\lambda_0 s}ds = \sum_{\pm} \frac{1- e^{-it(\lambda_0 \pm \sqrt{z})}}{2\sqrt{z}(\sqrt{z}\pm \lambda_0)}
\end{equation}
solves the corresponding ODE for each spectral component. Note that 
\begin{equation}\label{eq:wtzasymp}
    W_{t,\lambda_0}(z) \to (z - \lambda_0^2 + i0)^{-1} \quad \text{in } \mathcal{D}_z'(0, \infty)
\end{equation}
as $t\to \infty$. In order to apply \eqref{eq:wtzasymp} to $W_{t, \lambda_0}(z)$ we first write 
\begin{equation}
    W_{t, \lambda_0}(P) = \int_0^1 W_{t, \lambda_0}(z) dE(z)
\end{equation}
where $dE(z)$ is the spectral measure for $P$. Since \eqref{eq:wtzasymp} is a distributional limit, we need absolute continuity of $dE(z)$. By Stone's formula, the spectral measure for $P$ is given for $z\in [0,1]$ and $f\in C_c^\infty(\Omega)$ by 
\begin{equation}
    E(z)f = \frac{1}{2\pi i}\lim_{\epsilon\to 0^+} \Delta\left((P - (\sqrt{z}+i\epsilon)^2)^{-1} - (P - (\sqrt{z}-i\epsilon)^2)^{-1}\right)
\end{equation} 
Note that the proofs of theorem \ref{thm:lap} and theorem \ref{thm:diffinomg} hold if we take $\omega = \lambda -i\epsilon$, $\epsilon \downarrow 0$, in which case we get the outgoing solution $u(\lambda - i0)$. Thus theorems \ref{thm:lap} and \ref{thm:diffinomg} imply that
\begin{equation}
    E(z)f \in C^0((0,1); H^{1, -\alpha}_0(\Omega))\cap C^1((0,1); H^{0, -\beta}_0(\Omega))
\end{equation}
for any $\alpha>1/2, \beta>3/2$. In particular, the measure $\nu(z) = \Delta_{\Omega}^{-1}E(z)f$ is absolutely continuous with respect to the Lebesgue measure on $(0,1)$. Choose a cutoff function $\varphi\in C_c^\infty(\R; [0,1])$ such that $\varphi \equiv 1$ near $\lambda_0^2$ and $\text{supp}(\varphi)\subset (0,1)$. Then we split $W_{t,\lambda_0}(P)f$ into two terms:
\begin{align}
    w_1(t) &= \varphi(P)W_{t,\lambda_0}(P)f \\
    R_2(t) &= (1-\varphi(P))W_{t,\lambda_0}(P)f
\end{align}
We estimate $R_2(t)$ by the bound
\begin{flalign}
    \left|\Re(e^{i\lambda_0 t}(1-\varphi(z))W_{t, \lambda_0}(z))\right| &= |1-\varphi(z)|\frac{|\cos(t\sqrt{z}) - \cos(t\lambda_0)|}{|\lambda_0^2 - z|} \\
    &\leq \frac{2}{\lambda_0 \delta}
\end{flalign}
where 
\begin{equation}
    \delta = \text{dist}(\text{supp}(1-\varphi), \lambda_0^2)
\end{equation}
Then the functional calculus gives the uniform bound 
\begin{equation}
    \left\| \Re(e^{i\lambda_0 t}(1-\varphi(P))W_{t, \lambda_0}(P)f)\right\|_{H^1_0}\leq \frac{2}{\lambda_0 \delta}\|f\|_{H^{-1}_0}
\end{equation}
Next, we consider
\begin{equation}
    w_1(t) = \int_{\R}\varphi(z)W_{t, \lambda_0}(z)\nu(z)dz
\end{equation}
We now carefully show an analogue of the distributional limit \eqref{eq:wtzasymp} for $w_1(t)$. Using \eqref{eq:wtdef} we write $w_1(t)$ as
\begin{flalign}
    w_1(t) &= \frac{1}{2i}\int_0^t \int_{\R}\frac{\varphi(z)}{\sqrt{z}}\left(e^{is(\sqrt{z}-\lambda_0)} - e^{-is(\sqrt{z}+\lambda_0)}\right)\nu(z)dzds \\
    &= -i\int_0^t \int_{\R}\left(e^{is(\zeta - \lambda_0)} - e^{-is(\zeta + \lambda_0)}\right)H(\zeta)\varphi(\zeta^2)\nu(\zeta^2)d\zeta ds \\
    &= -i\int_{-t}^0\hat{q}_+(s)ds + i \int_0^t \hat{q}_-(s)ds
\end{flalign}
where 
\begin{equation}
    q_\pm(\zeta) = H(\zeta \pm \lambda_0)\varphi((\zeta \pm \lambda_0)^2)\nu((\zeta \pm \lambda_0)^2)
\end{equation}
and $H$ is the Heaviside function. Since $\nu$ is differentiable on the support of $\varphi$, $q_\pm\in H^1_c(\R, H^{0, -\beta}(\Omega))$. Therefore
\begin{flalign}
    \int_{\R}\|\hat{q}_\pm(s)\|_{H^{0, -\beta}(\Omega)} ds & \leq C\left(\int_\R \langle s \rangle^2 \| \hat{q}_\pm(s)\|_{H^{0, -\beta}(\Omega)}ds \right) \\
    & = \|q_\pm\|_{H^1(\R, H^{0, -\beta}(\Omega))}
\end{flalign} 
Thus, as $t\to \infty$, $w_1(t)$ converges in $H^{0, -\beta}(\Omega)$ to 
\begin{flalign}
    \lim_{t \to \infty}w_1(t) &= -i\int_{-\infty}^0\hat{q}_+(s)ds + i \int_0^\infty \hat{q}_-(s)ds \\
    &= \sum_{\pm} \int_\R(\zeta \pm i0)^{-1}q_{\pm}(\zeta)d\zeta \\
    &= \sum_{\pm}\int_\R (\zeta \pm (-\lambda_0 + i0))^{-1}H(\zeta)\varphi((\zeta + \lambda_0)^2)\nu((\zeta + \lambda_0)^2)d\zeta \\
    &= \sum_{\pm}\int_\R (\sqrt{z} \pm (-\lambda_0 + i0))^{-1}\frac{\varphi(z)}{\sqrt{z}}\nu(z)dz \\
    &= \int_{\R}(z-\lambda_0^2+i0)^{-1}\varphi(z)\nu(z)dz\\
    & = \varphi(P)(P-\lambda_0^2+i0)^{-1}f
\end{flalign}
So 
\begin{equation}
    u_1(t) := \Delta^{-1}w_1(t) = \Delta^{-1}(P-\lambda_0^2+i0)^{-1}f + \Delta^{-1}R_1(t) + \mathcal{E}(t)
\end{equation}
where 
\begin{equation}
    R_1(t) = (1-\varphi(P))(P-\lambda_0^2+i0)^{-1}f
\end{equation}
and $\mathcal{E}(t)\to 0$ in $H^{0, -\beta}(\Omega)$. Since $\varphi\equiv 1$ near $\lambda_0^2$, 
\begin{equation}
    \sup_{t}\|R_1(t)\|_{H^{-1}_0} \leq C\|f\|_{H^{-1}_0}
\end{equation}
with $C$ independent of $t$. Combining this bound with the estimate for $\Delta^{-1}R_2(t)$ gives the asymptotic \eqref{eq:timeasymp} for
\begin{equation}
    r(t)  = \Delta^{-1}R_1(t) + \Delta^{-1}R_2(t)
\end{equation}
and
\begin{equation}
u^+ = \Delta^{-1}(P-\lambda_0^2+i0)^{-1}f \in H_0^{1, -\alpha}(\Omega)
\end{equation}
which is the unique incoming solution to the stationary equation. \qed
\end{pf}

\nocite{ralston-1973,dyatlov-wang-zworski-2025,li2025evolutioninternalwaves2d}
\printbibliography
\end{document}